\documentclass[twoside,bezier]{amsart}

\usepackage{epsfig}
\usepackage{amsmath,amsthm,amsfonts,amssymb,amscd,euscript}

\input{epsf}
\newcommand{\filebegin}{\begin{document}}
\newcommand{\fileend}{\end{document}}
 \makeatother
\def\thefootnote{}
\newcommand{\lo}{\longrightarrow}
\newcommand{\NMM}{\hspace*{2mm}}

\newcommand{\cir}{\mathbb{T}}
\newcommand{\disk}{\mathbb{D}}

\renewcommand{\baselinestretch}{1.1}

\input{epsf}
\renewcommand{\baselinestretch}{1.1}
\def\n{\noindent}%

\numberwithin{equation}{section}
\def\mapdown#1{\Big\downarrow\rlap
{$\vcenter{\hbox{$\scriptstyle#1$}}$}}
\begin{document}
\vspace*{2cm}
\begin{center}

{\bf\large {New results on $p$-Carleson measures and some related
measures in the unit disk}}
 \\[0.5cm]
{Romi Shamoyan, Ali Abkar
\\ [0.5cm]
Department of Mathematics,
Bryansk State Technical University,
Bryansk 241050, Russia; email: rshamoyan@gmail.com \, and\\
Department of Mathematics, Imam Khomeini International University,
Qazvin 34149, Iran; email: abkar@ikiu.ac.ir\vspace{12pt}} \\[2mm]

\end{center}%
\vspace*{0.5cm}
\begin{quotation}
\noindent {\footnotesize {\sc Abstract. We provide some new sharp
embeddings for $p$-Carleson measures and some related measures in
the unit disk of the complex plane.\\}}
\end{quotation}

\markboth {R. Shamoyan, A. Abkar}
 {new results on $p$-Carleson measures and ...}
\section{Introduction}
The purpose of this  paper is to continue the study of $p$-Carleson
and related measures in the unit disk. In recent years many papers
appeared where $p$-Carleson measures were studied intensively from
various point of view by many authors (see [7], [8] and references
therein). We intend to study again $p$-Carleson measures
 in complex plane in simplest case, namely in the unit disk. Our
 intention in particular is to provide some direct generalizations
 of known one dimensional results about $p$-Carleson measures (or
 $Q_p$ analytic classes). On the other hand, the intention of this
 note is to provide sharp embeddings for certain new analytic spaces
 in the unit disk. This paper can be viewed as a continuation of
 [6] where similar assertions were recently proved.\\
Throughout the paper, we write $C$ (sometimes with indexes) to
denote a positive constant which might be different at each
occurrence (even in a chain of inequalities) but is independent of
the functions or variables being discussed. \\ Given two expressions
$A$ and $B$, we shall write $A\asymp B$ if there is a positive
constant $C$ such that $\frac B C\leq A\leq CB$. We will write
$A<_\prec B$ if there is a constant $C$ such that
$A<CB$. \\
In the next section we provide some new sharp embedding theorems for
so called $p$-Carleson measures and related measures in the unit
disk which are closely connected with $Q_p$ spaces in the unit disk
(see [7], [8]).\\
Let $\mathbb{D}$ be, as usual, the unit disk in the complex plane,
and let $\cir$ be the unit circle $\{ \xi :|\xi |=1\}$. By $dm_2$ we
mean the standard normalized Lebesgue area measure on $\disk$. We
denote the Lebesgue arc length measure on the unit circle by $d\xi$
or $dm(\xi)$. A positive Borel measure on $\disk$ will be denoted by
$\mu$. Let further $H(\disk)$ denote the space of holomorphic
functions in the unit disk. Moreover $$\square I=\{z=re^{i\varphi}:
1-|I|<r<1,\, \, \xi=e^{i\varphi}\in I\}$$ is a usual classical
Carleson box (see [6]) in $\disk$ where $I$ is an arc on $\cir$. Let
also $L^{p,q}(\cir),\, 0<p<\infty, 0<q<\infty$ be the classical
Lorentz space on $\cir$ (for definition see [1]) and
$$H^{q,\infty}(\cir)=H(\disk)\cap L^{q,\infty}(\cir).$$
By $\Gamma_\sigma(\xi)$ we mean the Lusin cone
$$\Gamma_\sigma(\xi)=\{ z\in\disk :|1-\overline \xi z|<\sigma(1-|z|)\},\quad
\sigma>1, \xi\in\cir.$$ Let
$$\triangle_z=\{z^{'}=r^{'}e^{i\theta^{'}}:|\theta -\theta^{'}|<1-r,\,
\frac{1-r}{2}<1-r^{'}<2(1-r)\},\quad z=re^{i\theta}.$$ For a fixed
$a\in\disk$, we let $\varphi_a$ be a M\"obius transformation on
$\disk$ interchanging zero and $a$ (that is,
$\varphi_a(z)=\frac{a-z}{1-\overline a z}$). The Bergman distance in
the unit disk which is M\"obius invariant is defined by
$$d(z,w)=\frac12\log \left (
\frac{1+|\varphi_w(z)|}{1-|\varphi_w(z)|}\right ),\qquad z,w\in\disk
.$$ Given a fixed $t>0$ and $a\in \disk$, the Bergman disk of radius
$t$ and center $a$ is defined to be (see [9]):
$$D(a,t)=\{ w\in
\disk :d(w,a)<t\}.$$ It is well known that the Bergman disk $D(a,t)$
is a Euclidean disk with center $\tilde a=\frac{1-r^2}{1-|a|^2r^2}a$
and the radius $R=\frac{1-|a|^2}{1-|a|^2r^2}r$ where
$r=\frac{e^{2t}-1}{e^{2t}+1}$; (see [9]). Note that
$$m_\alpha(D(a,t))\asymp(1-|a|)^{2+\alpha}$$ and
$dm_\alpha(z)=(1-|z|)^\alpha dm_2(z)$ for $\alpha >-1$ (for details,
see [9]). We also mention that for $a\in\disk$ and $z\in D(a,t)$ we
have
$$|1-\overline az|\asymp (1-|z|)\asymp (1-|a|).$$
It is well known that for every $\delta >0$, there exists a distinct
sequence $\{z_j\}$ in $\disk$ called a $\delta$-lattice, such that
$d(z_j,z_k)>\delta /5$ if $j\neq k$, and
$$\cup_jD(z_j,\delta)=\disk ;\qquad \sum_j \chi _{D(z_j,
5\delta)}(z) \le L,\quad z\in\disk,$$ where $L$ is a constant,
$\chi_E(z)$ is the characteristic function of a set $E$ (see [9]).

\section{Main Results} We now provide some new embedding theorems for $p$-Carleson
measures in the unit disk. In all our results, $\mu$ is a positive
Borel measure on $\disk$, and $f$ is an analytic function in $\disk$.\\
\noindent{\bf Theorem 1.} Let $q>1,\, \alpha >-1$, and $f$ be an
analytic function in $\disk$. Then the following conditions are
equivalent.\begin{enumerate}\item[(a)]
$$\left \| \int _{\Gamma _\sigma(\xi)}\frac{|f(z)|}{1-|z|}d\mu_ (z)\right \|
_{L^{q,\infty}(T)}\le C \left \| \int _{\Gamma
_\sigma(\xi)}|f(z)|(1-|z|)^{\alpha -1} dm_2 (z)\right \|
_{L^{q,\infty}(T)}.$$
\item[(b)] $\mu (D(z,t))\le C(1-|z|)^{\alpha +2},\quad t>0,\, z\in
\disk .$
\item[(c)] $\mu( \square I)\le C|I|^{\alpha +2}.$
\end{enumerate}
\begin{proof} We first show that (a) is equivalent to (b). For this
purpose we use the following two important known facts from [2] and
[4] (page 234):
\begin{equation}A(\mu)=\sup _{I\subset T} \frac{1}{|I|^s}\int
_{\square I}\frac{d\mu(z)}{1-|z|}\le C\sup_{I\subset
T}\frac{1}{|I|^s}\int_I\int_{\Gamma_\sigma(\xi)}\frac{d\mu(z)}
{(1-|z|)^2}=B(\mu)\end{equation} where $s\in (0,\infty)$ and
\begin{equation} A(\mu)\asymp B(\mu)\qquad \text{for}\quad
s\in(0,1),\end{equation}
\begin{equation}
\left(\|f\|^r_{L^{q,\infty}(T)}\right )= \sup _{I\subset T}
\frac{1}{|I|^{1-r/q}}\int _{I}|f(\xi)|^rdm(\xi),\qquad
0<r<q.\end{equation} As it was shown in [5] \begin{equation}
\frac{1}{(1-z\xi)^{1/q}}\in H^{q,\infty}(\cir),\quad \xi \in \cir ,
\end{equation}
and we also have (see [9])
\begin{equation}
 \int_{\Gamma_\sigma(\xi)}\frac{|f(z)d\mu(z)}{1-|z|}\le
\sup_{z\in \disk}\left(\frac{\mu(D(z,t))}{(1-|z|)^{\alpha +2}}\right
)\left(\int_{\Gamma_\sigma(\xi)}|f(z)|(1-|z|)^{\alpha
-1}dm_2(z)\right).
\end{equation}
Combining the relations (2.1) to (2.5) provided above, we obtain the
desired result. Indeed, one part follows directly from (2.5). For
the other part we note that our estimate in Theorem 1 by (2.2) and
(2.3) is equivalent to
\begin{multline*}\sup _{I\subset T}
\frac{1}{|I|^{1-1/q}}\int _{\square I}|f(z)|d\mu(z)\le\\
C\sup_{I\subset
T}\frac{1}{|I|^{1-1/q}}\int_I\int_{\Gamma_\sigma(\xi)}|f(z)|(1-|z|)^{\alpha
-1}dm_2(z).\end{multline*} Let
$$f(z)=\frac{(1-|w|)^\beta}{(1-zw)^{\beta +\alpha +1+1/q}},\qquad
z\in\disk,$$ where $\beta$ is a big enough positive number. Then the
left side is bigger that $\mu(\square I)/|I|^{\alpha +2}$ if $w$ is
a center of $\square I$, here $\square I$ is fixed and the right
side is bounded from above by some constant; here we can use the
fact that
$$|1-\overline wz|\asymp (1-|w|)\asymp |I|,\qquad z\in\square I,\,
w\in \disk,$$ where $w$ is a center of Carleson box, and
$$(1-|w|)^{-\beta}\le C|1-\overline wz|^{-\beta},\quad \beta >0.$$
Other implications are well known; see [7, 8, 9].
\end{proof}
We list two assertions concerning the actions of Lusin area integral
again, but we do not provide complete proofs for them, since the
main idea we have in this proofs coincides with the ideas mentioned
in the proof of Theorem 1. We have the following two statements.\\
\noindent{\bf Theorem 2.} Let $\beta>0, q>1, t>\beta +1$, and $f$ be
an analytic function in $\disk$. Then
$$\left \|\int _{\Gamma_\alpha(\xi)} \frac{(1-|z|)^t|f(z)|d\mu(z)}{1-|z|}\right \|
_{L^{q,\infty}(\cir)}\le C\sup_{|z|<1}|f(z)|(1-|z|^\beta)$$ if and
only if
$$\left \|\int _{\Gamma_\alpha(\xi)}(1-|z|)^{t-1}\left ( \int_{D(z,\tau)}d\mu(w)\right )
 \frac{(1-|z|)^{-\beta}}{(1-|z|)^2}\right \|
_{L^{q,\infty}(\cir)}\le C.$$
 \noindent{\bf Theorem 3.} Let
$1/p= 1/q+1/r,\, q, r, p >1,\, \beta >0, t>-1.$ Then for each $f\in
H(\disk)$ and each positive Borel measure $\mu$ on $\disk$ we have
\begin{multline*}\sup_{I\subset \cir}\frac{1}{|I|^{1-1/p}}\int_{\square
I}(1-|z|)^t|f(z)|d\mu(z)\le \\ C \sup_{I\subset
\cir}\frac{1}{|I|^{1-1/q}}\int_{ I}\sup _{z\in
\Gamma_\sigma(\xi)}|f(z)|(1-|z|)^\beta dm(\xi)\end{multline*} if and
only if
$$
\sup_{I\subset \cir}\left ( \int_{\square I}(1-|z|)^td\mu(z)\right
)\frac{1}{|I|^{\beta +1-1/r}} <\infty .$$ We now provide the
complete sketch of proofs for this last statements.\par {\bf Proof
of Theorem 2.} The idea of the proof is similar to the ideas
discussed in the proof of Theorem 1.
 One part of the theorem follows immediately from [9]:
\begin{multline*}\int _{\Gamma_\alpha(\xi)}(1-|z|)^{t-1}|f(z)|d\mu(z)\\ \le C \int
_{\Gamma_\alpha(\xi)}|f(z)|(1-|z|)^{t-1}\left
(\int_{D(z,\tau)}d\mu(w)\right )\frac{dm_\alpha (z)}{(1-|z|)^{\alpha
+2}},\end{multline*} for $\alpha >\beta -1,\, t>\beta +1.$\\
To show the reverse implication, we need to use the estimates
(2.1)-(2.4), the test function $f(z)=\frac{1}{(1-wz)^\beta}$ for a
large enough number $\beta$, and a fixed complex number $w$ in the
unit disk $\disk$. To be more precise, we use (2) to show that the
required estimate is equivalent to
$$\frac{1}{|I|^{1-1/q}}\int_{\square I}|f(z)|(1-|z|)^t d\mu(z)\le
C\sup_{|z|<1}|f(z)|(1-|z|)^\beta,$$ and also the facts that (see
[9]): \begin{multline*}\frac{C_1}{(1-|a|)^{\alpha +2}}\int_{D(a,
r+s)}d\mu(z)\\ \le \frac{1}{(1-|a|)^{\alpha
+2}}\int_{D(a,r)}(1-w)^{-\alpha -2}\int_{D(w,t)}d\mu(z)dm_\alpha(w)\\
\le \frac{C}{(1-|a|)^{\alpha +2}}\int_{D(a,r+t)}d\mu(z),\qquad t>s,
r>0.\end{multline*} This implies that the condition in the theorem
is equivalent to the fact that our measure is a $\tau$-Carleson
measure for some $\tau >0$. The rest of proof is routine.\qed \\
We now give a sketch for the proof of Theorem 3. First we verify the
necessity of the condition \begin{multline*}\sup_{I\subset
\cir}\frac{1}{|I|^{1-1/p}}\int_{\square
I}(1-|z|)^t|f(z)|d\mu(z)\le\\
C\sup_{I\subset \cir} \frac{1}{|I|^{1-1/q}}\int_I
\sup_{z\in\Gamma_\sigma(\xi)}|f(z)|(1-|z|)^\beta
dm(\xi).\end{multline*} Since $(1-z)^{-1/q}\in H^{ q,\infty}$ (see
[5]) we can put $f(z)=\frac{1}{(1-z)^{\beta +1/q}}$. The right side
is bounded. As for the left hand side of the estimate above, we can
estimate it from below. Indeed, we have to verify that the following
estimate is true. We have
\begin{multline*}\sup_{I\subset
\cir}\frac{1}{|I|^{1-1/p}}\int_{\square
I}\frac{(1-|z|)^td\mu(z)}{|1-\overline{w} z|^{\beta+1/q}} \ge \\
C\sup_{I\subset \cir}\left( \int_{\square I}(1-|z|)^td\mu(z)\right
)\frac{|I|^{-\beta}}{|I|^{1-1/r}}\end{multline*} where $w$ is a
center of $\square I$, and $1/p=1/q +1/r,\, q, r,  p
>1,\, \beta >0,\, t>-1.$ This is true, since if $w$ is a center of
$\square I$, then
$$|1-\overline wz|\asymp (1-|w|)\asymp |I|,\qquad z\in\square I,\,
w\in \disk.$$ Hence we obtain
$$(1-|w|)^{-\beta}\le C|1-\overline wz|^{-\beta},\quad \beta >0.$$
To show the reverse, we have to modify the proof of Theorems 1,2 and
use the estimates (2.1)-(2.5). We have
\begin{multline}
\int_{\Gamma _\sigma(\xi)}\frac{(1-|z|)^t|f(z)|d\mu(z)}{1-|z|}\le
C\left (
\sup_{z\in\Gamma_\sigma(\xi)}|f(z)|(1-|z|)^\beta \right )\times \\
\times \int_{\Gamma_\sigma(\xi)}\frac{(1-|z|)^t}{(1-|z|)^{\alpha
+2}}\left (\int_{D(z,r)}d\mu(w) \right )\frac{dm_{\alpha
-\beta}(z)}{1-|z|},\qquad \beta >0, t>0.\end{multline} Then we use
(see [4])
$$\| f_1f_2\|_{L^{ p,\infty}(\cir)}\le C
\| f_1\|_{L^{ q,\infty}(\cir)} \| f_2\|_{L^{ r,\infty}(\cir)},\quad
\frac {1}{ p}= \frac {1}{ q}+ \frac {1}{ r},\,  q,  r,  p >1.$$
Using (2.2) in (2.6) we will have the estimate
\begin{multline*}\sup_{I\subset
\cir}\frac{1}{|I|^{1-1/p}}\int_{\square
I}(1-|z|)^t|f(z)|d\mu(z)\le\\
C\sup_{I\subset \cir} \frac{1}{|I|^{1-1/q}}\int_I
\sup_{z\in\Gamma_\sigma(\xi)}|f(z)|(1-|z|)^\beta dm(\xi)\times
\\\sup_{I\subset\cir}\frac{1}{|I|^{1-1/r}}\int_{I}
\int_{\Gamma_\sigma(\xi)}\frac{(1-|z|)^t}{(1-|z|)^{\alpha +2}}\left
(\int_{D(z,r)}d\mu(w) \right )\frac{dm_{\alpha
-\beta}(z)}{1-|z|}=G(f)M(\mu).
\end{multline*}
It can be easily shown that \begin{multline*}M(\mu)\le C
\sup_{I\subset\cir}\frac{1}{|I|^{1-1/r}}\int_{\square
I}(1-|z|)^{t-\alpha -2}\left (\int_{D(z,r)}d\mu(w) \right
)dm_{\alpha -\beta}(z)\\ \le C  \sup_{\tilde z\in
\disk}\frac{1}{(1-|\tilde z|)^{1-1/r}}\int_{D(\tilde z,
r_1)}(1-|z|)^{t-\alpha -2}\left (\int_{D(z,r)}d\mu(w) \right
)dm_{\alpha -\beta}(z).\end{multline*} Put $1-1/r =\alpha +2$, then
it follows that $\alpha =-1-1/r$. Therefore
\begin{multline*}M(\mu)\le
C  \sup_{\tilde z\in \disk}\frac{1}{(1-|\tilde z|)^{\alpha
+2}}\int_{D(\tilde z, r_1)}\frac{1}{(1-|z|)^{\alpha +2}}\\
\times \left (\int_{D(z,r)}(1-|\tilde z|)^{t-\beta}d\mu(\tilde z)
\right )dm_\alpha(z)\\ \le C\sup \left ( \int_{D(\overline z,
r_1+r)}(1-|\tilde z|)^td\mu (\tilde z_1)\right )\frac{1}{(1-|\tilde
z_1|)^{\alpha +2+\beta}}=M_1(\mu).\end{multline*}
 \noindent{\bf Theorem 4.} Let $0<p<q, \, \alpha >0$, and $f$ be an analytic
function in $\disk$. Then
$$\sup_{I\subset\cir}\frac{1}{|I|^{1-p/q}}\int_{I}\int_{\Gamma_\sigma
(\xi)}\frac{|f(z)|^pd\mu(z)}{(1-|z|)}\le C \left
(\int_{\disk}|f(z)|^q(1-|z|)^{\alpha q-1}dm_2(z) \right )^{1/q}$$ if
and only if
$$\sup_{I\subset\cir}\frac{1}{|I|}\int_{\square I}\left (\frac{\mu
(\Delta _z)}{(1-|z|)^{\alpha p+1}} \right
)^{\frac{q}{q-p}}\frac{dm_2(z)}{1-|z|}\le C.$$ Ideas for the proof
of this theorem are the same as those used in the proof of Theorem
2, so that we omit the details.\par We remark that the conditions on
measures appeared in the formulations of Theorems 2 and 3 can be
reformulated in terms of some conditions on $D(z, t)$, as we have
done in Theorem 1.


\providecommand{\bysame}{\leavevmode\hbox
to3em{\hrulefill}\thinspace}

\end{document}